\documentclass[pdflatex,sn-basic]{sn-jnl}
\setcitestyle{aysep={,}}

\usepackage{
    amsmath,
    amssymb,
    commath,
    dsfont,
    cleveref,
    booktabs,
    multirow,
    algorithm,
    algpseudocodex,
}

\usepackage[
    group-separator={,}
]{siunitx}

\DeclareMathOperator*{\argmin}{argmin}
\DeclareMathOperator*{\argmax}{argmax}

\renewcommand{\P}{\mathds{P}}
\newcommand{\diff}{\mathop{}\!\mathrm{d}}

\title{Simulation Strategies for an Efficient Local Search to solve Stochastic Scheduling Problems}

\author*[1]{\fnm{Philip} \sur{de Bruin}}\email{p.debruin@uu.nl}
\author[1]{\fnm{Bram} \sur{Elderhorst}}\nomail
\author[1]{\fnm{Marjan} \sur{van den Akker}}\email{j.m.vandenakker@uu.nl}
\author[1]{\fnm{Han} \sur{Hoogeveen}}\email{j.a.hoogeveen@uu.nl}

\affil[1]{\orgdiv{Department of Information and Computing Sciences}, \orgname{Utrecht University}, \country{The Netherlands}}

\abstract{%
In scheduling problems, deterministic task durations are often assumed.
This usually does not capture reality and may lead to schedules that are not robust to (small) changes to these task lengths.
The use of stochastic task durations therefore seems preferable.
Including these in local search, which is the way to find good solutions for difficult scheduling problems, is not straightforward, though.
The objective value becomes stochastic then too, and computing the expected value is often not possible.
One way out of this it to approximate this value by using simulation.
This is quite easy to implement in a local search algorithm, but it may require many simulations each iteration to get a reliable estimate.
Hence such an approach comes with a performance penalty.

In this paper, we study techniques to limit the number of simulations.
Besides comparing known techniques, we propose our own method for this, which is based on $t$-tests.
We evaluate these techniques on the \textsc{Stochastic Parallel Machine Scheduling Problem} and the \textsc{Stochastic Electric Vehicle Scheduling Problem}.
In these case studies, we show the effectiveness of using such methods to reduce runtime while retaining solution quality.
Our method using $t$-tests turns out to be most effective in both problems.
}


\begin{document}

\maketitle

\section{Introduction}\label{sec:introduction}
In scheduling, it is quite common to assume the length of tasks to be deterministic.
While many solution algorithms are known for deterministic scheduling problems, it does not capture reality really well.
For example, in bus scheduling, one can assume the driving times to be deterministic, but in reality buses have to deal with different traffic conditions causing small deviations.
Thus, the use of deterministic driving times may result in non-robust schedules.
This comes at the cost of worse on-time performance and passenger satisfaction, showing the need to better deal with these variations.
To achieve this, driving times can be modelled as stochastic variables.
Difficult scheduling problems can be successfully solved using heuristics such as local search.
In this paper, we focus on how to incorporate stochastic variables into local search.

The inclusion of stochastic variables in a local search algorithm can be done in several ways.
Fully incorporating these stochastic variables, i.e. formulating the objective function in terms of probability distributions, can prove difficult if it is very hard to compute the resulting distribution.
Another possibility is to convert the stochastic variables to deterministic values.
For example, one could use a certain percentile of the underlying distribution, or multiply the mean by a certain factor, which we call the robustness factor.
This approach has the benefit that one could still use deterministic values when optimizing, however this does not fully account for the variance in the distribution.
It especially does not account for bigger disturbances that realistically could still occur.

Instead of converting the stochastic variables to deterministic values, one could also create a deterministic function to approximate the objective function of a given solution.
An example of these are so-called surrogate robustness measures, which are shown to be effective in scheduling problems by \citet{Loman2025}.

Another way to approximate the objective value is by using simulation.
The inclusion of simulation in an optimization algorithm is generally known as simulation-based optimization.
We specifically focus on including simulation in the Simulated Annealing algorithm, however most techniques can be applied in other local search settings as well.
\Citet{Alkhamis1999} and \citet{Alrefaei1999} were one of the first to provide frameworks for using simulations inside a Simulated Annealing algorithm.

Both \citeauthor{Alkhamis1999} and \citeauthor{Alrefaei1999} note that the use of simulations to evaluate the objective function may have a big impact on runtime performance.
In order to get a good estimate of the quality of a solution, a lot of simulations may have to be performed.
This is especially true when there is a larger variance between simulations.
One way to combat this is by employing techniques from \emph{Ordinal Optimization} \citep{Ho1999}.
In ordinal optimization one has to create a correct ranking of different solutions.
This is similar to what is happening in local search algorithms that look for iterative improvements.
These local search algorithms compare in each iteration two solutions to determine which one is better.
Hence, it is not necessary to know the exact score of a solution.
It is shown that using these ordinal optimization techniques can greatly improve the performance of various heuristic approaches \citep{Lee2010,Yang2014}.
These speed-ups are generally realized by not spreading the simulations evenly over all solutions, but to use more simulations on some solutions and fewer on others.
Specific techniques that do this are \emph{Optimal Computing Budget Allocation} (OCBA) and so-called \emph{Indifference Zones} (IZ) \citep{Chen2000,Kim2006}.

These ordinal optimization methods are quite generic in the sense that they compare $k$ solutions.
Furthermore, they are designed to create a ranking between these solutions, while local search algorithms such as Simulated Annealing also take the difference between solutions into account, as this information is used to determine the probability of selecting a worse solution to escape local optima.
Our contribution is twofold.
First, we propose a new method based on $t$-tests to calculate the required number of simulations with the purpose to reduce the total number of simulations.
This method can create a ranking between two solutions, and it can better incorporate the difference between solutions.
We include this new method and other ordinal optimization methods, namely OCBA and our adaptation of IZ, in a Simulated Annealing algorithm.
Second, we present two case studies to evaluate and compare the resulting algorithms.
The first case study focuses on the \textsc{Stochastic Parallel Machine Scheduling Problem}, where jobs have release dates, precedence relations, and stochastic processing times.
The second case study is on the \textsc{Electric Vehicle Scheduling Problem} for public bus transport.
Here, both stochastic driving times and stochastic energy consumption are modelled using historic data.

The paper is structured as follows.
In \Cref{sec:literature} we present a literature overview.
Then in \Cref{sec:simulated_annealing} we show how the Simulated Annealing is adapted to incorporate these stochastic variables.
We show the various approaches to determine the number of simulations in each iteration in \Cref{sec:number_simulations}.
We provide computational results for each case study in \Cref{sec:cs_spmsp,sec:cs_evsp}, and the conclusion in \Cref{sec:conclusion}.

\section{Literature Overview}\label{sec:literature}
Stochastic optimization problems are optimization problems with a stochastic objective and/or constraints.
There are several ways to solve these problems.
One method of solving is using simulation-based optimization, where simulation is used to estimate the stochastic parts of the model.
\Citet{Amaran2016} provide a survey on the field of simulation-based optimization.
Next to that, they show the relationship of this field to other well-known methods, such as robust optimization \citep{Bertsimas2011} and stochastic programming \citep{Birge1997}.
\Citet{Amaran2016} describe several techniques within simulation-based optimization.
Examples of explored techniques are: ordinal optimization, random search methods, and direct search methods.
In ordinal optimization, the focus lies on sampling from a subset of solutions and finding the best among them.
The key idea here is to focus on the ranking between these solutions instead of having a good estimate of the objective value.
Direct search methods iteratively search for a better solution, where the iterations have a deterministic pre-specified pattern \citep{Kolda2003}.
Random search methods include metaheuristics such as Simulated Annealing and Tabu Search.
A survey specifically on the use of metaheuristics for stochastic optimization problems is presented by \citet{Bianchi2009}.
They describe two ways to integrate metaheuristics into stochastic optimization problems.
Either by designing an approximation function for the objective, or by simulation.
More recently, the term `simheuristics' was introduced to describe the integration of simulation and metaheuristics \citep{Juan2022}.
An overview of the use of simheuristics in different problem settings is provided by \citet{Juan2015}.

\Citet{Alkhamis1999} and \citet{Alrefaei1999} are the first to provide a framework in which they combine Simulated Annealing with an objective value that is estimated from many samples, i.e. simulations.
Both papers use a given number of $N_k$ simulations in iteration $k$.
\Citet{Alkhamis1999} prove convergence conditions for this approach, showing that essentially the standard deviation of the sampled objective value has to converge to zero as $k \to \infty$.
They note that this criterion might not be very practical, as it may require many simulations to sufficiently reduce the standard deviation.
\Citet{Alrefaei1999} provide computational experiments on optimization problems in an $M / M / 1$ queuing system.
They note that their approach yields good performance when they use a slowly increasing sequence $\{ N_k \}$ to estimate the objective value.
These approaches were later extended by \citet{Alkhamis2004}, who incorporate the standard deviation of the objective value into the acceptance criterion of Simulated Annealing.

\Citet{Bianchi2009} remark that while the use of Simulated Annealing for stochastic optimization has received quite some attention, the research mainly focuses on understanding the conditions for convergence to the optimal solution.
In that regard, not much is done when it comes to applying these techniques to actual problems.
\Citet{Akker2013} use such an approach to find robust solutions for the \textsc{Stochastic Job Shop Scheduling Problem}.
They keep the number of simulations in each iteration constant during the whole search process.
A similar approach was used by \citet{Tasoglu2019} to solve the integrated berth allocation and quay crane scheduling problems with stochastic handling times.

\Citet{Passage2025} investigate the \textsc{Stochastic Parallel Machine Scheduling Problem}.
They use a Variable Neighbourhood Descent method in combination with Iterated Local Search.
Instead of using simulations to sample their objective value, they estimate it by assuming every stochastic variable to be normally distributed, which makes it possible to compute the expected value and variance of each variable, including the one corresponding to the objective value.
Their approach is compared to the simulation sampling approach of \citet{Akker2013}.
They find that computing these estimates is very efficient, especially when compared to using simulation.
Using \num{300} to \num{1000} simulations per iteration outperforms their approach, although with larger computation times.
Here, the exact number of simulations needed depends on the distribution of the variables.

One way to improve the time spent on simulations in a local search algorithm is by minimizing the number of simulations needed while still doing enough simulations to make correct decisions.
In a local search framework, in each iteration different solutions are compared to each other in order to select the `right' one.
Local search algorithms that look for iterative improvement always select the best solution, while Simulated Annealing can select a worse solution.
Thus, in case of iterative improvement we would mostly be interested in the correct ordering of solutions instead of an accurate score of each solution.
In the Simulated Annealing setting, however, we need to know the score difference between the two solutions to determine whether the neighbour is accepted, since the probability of accepting a worse solution depends on this.
Thus knowing a correct ordering is not sufficient.
Still techniques for finding a correct ordering could serve as a basis for reducing the number of simulations in a Simulated Annealing setting.
Finding a correct ordering of $k$ solutions is called \emph{Ordinal Optimization} \citep{Ho1999}.
One technique here is \emph{Optimal Computing Budget Allocation} (OCBA), which maximizes the probability of a correct selection given a maximum number of simulations.
This method works by allocating simulations to a solution based on the current sample variance instead of dividing the simulations evenly.
OCBA has been used by \citet{Yang2014} in the \textsc{Stochastic job Shop Scheduling} problem and by \citet{Clapper2024} in stochastic home healthcare routing and scheduling.
Another technique for \emph{Ordinal Optimization} is by using \emph{Indifference Zones} \citep{Lee2010, Ahmed2002, Ghiani2007}.
This is similar to OCBA, but where OCBA maximizes the probability of a correct selection, Indifference Zones guarantee a minimum for the probability of a correct selection \citep{Lee2010}.

\section{Simulated Annealing Framework}\label{sec:simulated_annealing}
In this section, we give an overview of how the Simulated Annealing algorithm is adapted to use simulations for the objective function.
Here, we look at a minimization problem.
In a standard Simulated Annealing implementation, we take an initial solution $s_\text{init}$ and run a number of iterations until we reach some stopping criterion.
Each iteration, a neighbour $s_n$ is generated based on the current solution $s$.
The Boltzmann equation is used to decide whether to accept or decline this neighbouring solution.
For this, the temperature parameter $T$ is used, which is decreased every $Q$ iterations.
In this minimization problem, given a solution $s$ with cost $c_s$, a neighbouring solution $s_n$ with cost $c_{s_n}$, and a temperature $T$, the neighbouring solution is accepted with probability
\begin{equation}\label{eq:sa_acceptance}
    \P(c_s, c_{s_n}, T) = \begin{cases}
        1 & \text{if $c_s \geq c_{s_n}$}, \\
        \exp\left( \frac{c_s - c_{s_n}}{T} \right) & \text{otherwise}.
    \end{cases}
\end{equation}

With this probability scheme, the Simulated Annealing does not only accept better solutions, but can accept a worse solution.
This can be an issue when determining a relative order of the solutions.
To resolve this, we can calculate how much worse the new solution is allowed to be at the start of an iteration.
This information could then be used when determining the relative ordering.
More specifically, we sample a random number $u \in [0, 1]$ at the start of an iteration.
The new solution $s_n$ is accepted if $u \leq \exp\left( \frac{c_s - c_{s_n}}{T} \right)$.
Note that this is true for both cases of \Cref{eq:sa_acceptance}, since $\exp\left( \frac{c_s - c_{s_n}}{T} \right) \geq 1$ if $c_s \geq c_{s_n}$.
This can be rewritten to $c_{s_n} + T\ln(u) \leq c_{s}$.
We call $T\ln(u)$ the allowed difference and denote this by $D$.
Note that $D < 0$.

Recall that in a stochastic problem setting, it is often not feasible or possible to calculate the cost $c_s$.
We combat this by estimating $c_s$ using simulations.
For this estimate we take the average cost of the simulations, which we denote by $\hat{c}$.
Thus, to accept or decline a neighbouring solution, we need to compare $\hat{c}_s$ and $\hat{c}_{s_n}$.
The pseudocode of the Simulated Annealing is given in \Cref{alg:simulated_annealing_with_simulation_estimation}.
In this algorithm, there are two moments where we need to determine how many simulations are performed, namely when comparing the new neighbour against the current solution and when comparing a newly accepted solution against the best one found so far.
For these cases we use $N$ and $N_b$ simulations respectively.
Note that in our setup the values for $N$ and $N_b$ are determined using the same method.
It is also possible to use different methods here, since it might be more important to make a correct decision on updating $s_\text{best}$.
However, our preliminary experiments found that this is not necessary and hence we use the same method for determining both $N$ and $N_b$.

The parameters used in the Simulated Annealing are the initial temperature $T_\text{init}$, the cooling rate $\alpha_\text{cool}$, and the number of iterations per cooling step $Q$.
We use an exponential cooling scheme, where every $Q$ iterations the current temperature $T$ is multiplied by $\alpha_\text{cool}$.
Furthermore, in both our case studies the Simulated Annealing is stopped once a certain temperature is reached.

\begin{algorithm}
    \caption{Simulated Annealing algorithm with simulations}
    \label{alg:simulated_annealing_with_simulation_estimation}
    \begin{algorithmic}[1]
        \State $T \gets T_{\text{init}}$
        \State $s \gets s_{\text{init}}$
        \State $s_\text{best} \gets s$
        \While{stop condition is not met}
            \If{current iteration number is divisible by $Q$}
                \State $T \gets \alpha_\text{cool} \cdot T$
                \Comment{Update temperature}
            \EndIf
            \State{\color{white}.} 
            
            \State $s_n \gets$ \Call{Neighbour}{$s$}
            \Comment{Generate a neighbouring solution}
            \State $u \gets$ \Call{Random}{0,1}
            \State $D \gets T\ln(u)$
            \Comment{Determine allowed difference}
            \State{\color{white}.} 

            \State calculate $\hat{c}_s$ and $\hat{c}_{s_n}$ using $N$ simulations, with $N$ prescribed by the considered method \label{line:estimate_cost_determine_n}
            
            \If{$\hat{c}_{s_n} + D \leq \hat{c}_s$}
                \State $s \gets s_n$
                \Comment{Accept $s_n$}
                \State{\color{white}.} 
                
                \LComment{Check if best solution needs to be updated}
                \State calculate $\hat{c}_s$ and $\hat{c}_{s_\text{best}}$ using $N_b$ simulations, with $N_b$ prescribed by the same method as on Line \ref{line:estimate_cost_determine_n} \label{line:estimate_cost_best_solution}
                \If{$\hat{c}_s \leq \hat{c}_{s_\text{best}}$}
                    \State $s_\text{best} \gets s$
                    \Comment{Update best solution}
                \EndIf
            \EndIf
        \EndWhile
    \end{algorithmic}
\end{algorithm}

To create a `fair' comparison between solutions, one can employ a technique called \emph{Common Random Numbers} (CRN) \citep{Law2015}.
With this, we make sure that both solutions are simulated on the same realizations of random numbers, thus one solution cannot get `lucky' by drawing `better' random numbers.
Note, however, that the use of CRN might not always be possible.
This is the case when there are stochastic variables that depend on the solution structure.
For example, \citet{Bruin2023} model stochastic energy consumption of electric buses, where they use a high-level model to determine which driver is driving the bus in order to have more realistic energy consumptions patterns.
Their approach means that the randomly chosen energy consumption of a trip depends on the previously driven trips.
Since this trip sequence will be different between solutions, the use of CRN for these variables is not possible in such cases.
Note that in these cases it is still possible to use CRN for the variables that do not depend on the solution structure.

A drawback of CRN is that it impacts the performance when updating the best found solution, i.e. we cannot re-use the result of $\hat{c}_{s_\text{best}}$ on Line \ref{line:estimate_cost_best_solution} and need to reevaluate the cost.
When calculating $\hat{c}_{s_\text{best}}$ on Line \ref{line:estimate_cost_best_solution}, we can only use a previous result if we evaluate on the exact same scenarios (samples) as in previous iterations.
\Citet{Passage2025} show that using the same set of scenarios helps with convergence, however this means that one needs to select enough scenarios with enough variation in order to not get solutions that only perform well on the selected scenarios.
Next to that, we compare methods that do not perform the same number of simulations each iteration, hence having a fixed set of scenarios is not suited for these methods.
Hence, we will be using a different set of scenarios each iteration.

\section{Number of Simulations}\label{sec:number_simulations}
In this section, we describe different methods to determine the number of simulations $N$ (Line \ref{line:estimate_cost_determine_n} of \Cref{alg:simulated_annealing_with_simulation_estimation}).
Here, it is important to have an accurate cost estimate, since otherwise the Simulated Annealing will have difficulties converging to an optimal solution.
More simulations, i.e. higher $N$, results in a more accurate estimate, but comes with a performance penalty.
Note here that two solutions are easier to tell apart when their score difference is bigger.
In these cases, fewer simulations will be required.
Thus, the goal becomes to minimize the number of simulations while still performing enough simulations such that the Simulated Annealing algorithm can make a decision with enough (the desired level of) accuracy.

We first look into Optimal Computation Budget Allocation, which is a method to divide a given number of simulations over the solutions such that the reliability of the decision is maximized.
Next, we look into two methods that determine the number of simulations dynamically, instead of using a constant number of simulations each iteration.
This should create a good balance between simulating too much, sacrificing computation time, and simulating too little, potentially sacrificing overall solution quality.
These methods are Indifference Zoning, and our new technique using paired $t$-tests.
In the remainder of this section, these three techniques are explained in further detail, including how they are adapted to work within Simulated Annealing algorithm, specifically how the method can be adapted to use information about the allowed difference $D$.
Note that all these methods just decide when to stop simulating.
The decision of which solution is accepted is still handled by the Simulated Annealing algorithm and is based on the estimated averages.

\subsection{Optimal Computation Budget Allocation}\label{sec:oo_ocba}
The general idea of Optimal Computation Budget Allocation (OCBA) is to divide a budget of simulations between the solutions we compare such that we maximize the probability of a correct selection.
Thus, this method always uses the same number of simulations, but instead of dividing these simulations evenly over all solutions, it divides them such that the probability of correct selection is maximized.
We implemented the algorithm described by \citet{Chen2000}.
This algorithm defines three parameters: $n_0$, $\Delta$, and $N_\text{max}$.
These are the initial number of simulations, the budget increase per step, and the maximum budget, respectively.
The algorithm first simulates each solution $n_0$ times.
It then repeatedly distributes $\Delta$ simulations over the solutions until we hit the computation budget $N_\text{max}$.

The additional $\Delta$ simulations are allocated using the following rule.
First note that in the Simulated Annealing, we are comparing two solutions with each other.
Let $N_1$ be the number of simulations we will use for the first solution and $N_2$ the number of simulations we will use for the second solution.
Thus, after the initial simulations, we have $n_0 = N_1 = N_2$.
Furthermore, let $s^2_1$ and $s^2_2$ denote the sample variances for each of the corresponding solutions.
Without loss of generality, we assume that solution $1$ is currently better than solution $2$.
\Citet{Chen2000} find that the optimal allocation of simulations should follow
\begin{equation}
    \frac{N_1}{N_2} = \frac{s_1}{s_2}
.\end{equation}
Hence, we allocate the additional $\Delta$ simulations in such a way that $\frac{N_1}{N_2}$ better approximates the ratio $\frac{s_1}{s_2}$.
Or, in other words, we minimize the difference between $\frac{N_1 + \Delta - i}{N_2 + i}$ and $\frac{s_1}{s_2}$, where $0 \leq i \leq \Delta$ and $i$ integer.

A potential drawback of this allocation framework is that it does not allow us to implement CRN.
For CRN we would need both solutions to get the same number of simulations, which this framework does not guarantee.

Note that this method only uses the sample variance information to divide the simulations.
This means that this method cannot be adapted to use information about the allowed difference.
However, in Simulated Annealing decisions are made based on the sampled averages and hence it can still look at the difference between solutions.
This might result in less accurate decisions being made by the simulated annealing, but it remains a question if that is a big issue.

\subsection{Indifference Zones}\label{sec:oo_iz}
Indifference Zones (IZ) work differently compared to OCBA.
They guarantee a minimum probability of correct selection, instead of maximizing this probability.
A drawback of this is that we do not have a maximum number of simulations per iteration, as we cannot guarantee that the required minimum probability of correct selection is reached within a certain number of simulations.
We look at a minimization problem and define the width of the indifference zone $\delta^*$ and a desired confidence value $\alpha_\text{conf}$.
If two solutions are within $\delta^*$ units of each other, the decision maker considers them to be the same, or ``indifferent''.
Then IZ procedures guarantee the following \citep{Kim2006}:
\begin{equation}\label{eq:iz_guarantee}
    \P(\text{CS}) = \P(\text{Solution 1 is observed as best} \mid \mu_1 + \delta^* \leq \mu_2) \geq 1 - \alpha_\text{conf}
.\end{equation}
Here, $\P(\text{CS})$ denotes the probability of correct selection, and $\mu_1$ and $\mu_2$ denote the (unknown) true means of solutions 1 and 2 respectively.

One way to guarantee this probability of correct selection is by following Rinott's two-stage procedure \citep{Rinott1978}.
This procedure compares $k$ different solutions and determines which one is best.
It first simulates each solution $n_0$ times and calculates the sample variances of each solution $s^2_i$.
With this information, the number of simulations needed for solution $i$ is calculated as follows:
\begin{equation}
    N_i = \max\left\{
        n_0,
        \left\lceil \left( \frac{h s_i}{\delta^*} \right)^2 \right\rceil
    \right\}
.\end{equation}
Here, $h$ denotes the solution to Rinott's double integral equation
\begin{equation}
    \int_0^\infty \left[ \int_0^\infty \Phi\left( \frac{h}{\sqrt{(n_0 - 1) (\frac{1}{x} + \frac{1}{y})}} \right) f_{n_0 - 1}(x) \diff x \right]^{k - 1} f_{n_0 - 1}(y) \diff y
    = 1 - \alpha_\text{conf}
,\end{equation}
where $\Phi$ denotes the cumulative probability function of the standard normal distribution, and $f_n$ denotes the probability density function of the $\chi^2$ distribution with $n$ degrees of freedom.
On modern computers, this integral can be solved quickly, but tables for $h$ are also available in \citep{Wilcox1984}.
Lastly, each solution gets simulated an additional $N_i - n_0$ number of times, after which the solution with the lowest mean is selected as best.

\Citet{Yoon2019} remark that Rinott's procedure does not use sample mean information, resulting in conservative results for the total number of simulations.
For this reason, they change the procedure into an iterative approach where in each iteration they use the sample means to decide which subset of the $k$ solutions requires additional simulations.
This is done until no solution needs more simulations.
Their procedure is as follows:
\begin{enumerate}
    \item Simulate each of the $k$ solutions $n_0$ times and calculate the sample means $\bar{x}_i$ and sample variances $s^2_i$.
    Furthermore, let $I$ be the set containing all $k$ solutions, let $N_i = n_0$, and $b = \argmin_i \bar{x}_i$.

    \item\label{item:my_stopping_condition}
    The set $I$ will now be updated such that every solution remaining in $I$ after this step gets additional simulations.
    For this, calculate $\delta_i = \max\{ \delta^*, \bar{x}_i - \bar{x}_b \}$, remember that $\delta^*$ is the zone width in which two solutions are `indifferent' to each other.
    Furthermore, $h_1$ is the solution to the integral
    \begin{equation}
        \int_0^\infty \left[
            \int_0^\infty \Phi\left(
                \frac{h_1}{\sqrt{(N_i - 1)\frac{1}{x} + (N_b - 1)\frac{1}{y}}}
            \right) f_{N_i - 1}(x) \diff x
        \right] f_{N_b - 1}(y) \diff y
        = 1 - \frac{\alpha_\text{conf}}{k - 1}
    ,\end{equation}
    where $f_n(x)$ denotes the probability density function of the $\chi^2$ distribution with $n$ degrees of freedom.
    
    Then, delete solution $i$ ($i \not= b$) from $I$ if
    \begin{equation}
        N_i \geq \left\lceil \left( \frac{h_1 s_i}{\delta_i} \right)^2 \right\rceil
        \text{ and }
        N_b \geq \left\lceil \left( \frac{h_1 s_b}{\delta_i} \right)^2 \right\rceil
    ,\end{equation}
    and delete solution $b$ from $I$ if
    \begin{equation}
        N_b \geq \left\lceil \left( \frac{h_1 s_b}{\delta_i} \right)^2 \right\rceil
        \text{ for all $i \not= b$}
    .\end{equation}

    \item If $\abs{I} = 0$, stop and return solution $b$ as the best solution.

    \item\label{item:my_more_simulation}
    Else, give each solution in $i \in I$ one additional simulation and set $N_i \gets N_i + 1$.
    Then, restore $I$ to contain all $k$ solutions again, and update the sample means $\bar{x}_i$ and sample variances $s_i^2$.
    Also update $b = \argmin_i \bar{x}_i$ and go back to Step \ref{item:my_stopping_condition}.
\end{enumerate}
\Citet{Yoon2019} show that \Cref{eq:iz_guarantee} holds for the result of this procedure, while needing significantly fewer simulations compared to Rinott's procedure.
Note that in this algorithm, the value of $h_1$ can be pre-computed (as can the value $h$ in Rinott's double integral).

We adapted this procedure to improve the performance when implemented within Simulated Annealing.
In Step \ref{item:my_more_simulation}, instead of giving each solution in $I$ one additional simulation, we give each solution $\Delta$ additional simulations, with $\Delta$ a parameter to choose.
While this could lead to more simulations than needed, overall it improved the runtime performance significantly.

We further improved the runtime performance by introducing a maximum number of simulations to be run per local search iteration.
This does break the guarantee in \Cref{eq:iz_guarantee}.
In preliminary experiments, however, we did not see a significant score difference between the case with and without this maximum.
This is likely due to the fact that Simulated Annealing can choose a worse solution with a certain probability, meaning that in cases where the solution costs are close to each other, the probability of correct selection might have a minor impact.

Lastly, the procedure does not allow for CRN, since there is no guarantee that each solution gets an equal number of simulations.
This can be solved by giving all solutions, rather than just the ones in $I$, additional simulations in Step \ref{item:my_more_simulation}.
Since we are only comparing two solutions ($k = 2$), we expect that this does not impact the runtime performance too much, while we get the benefits of having CRN.

To use the allowed difference information in the Simulated Annealing, we compare $\bar{x}_{s}$ with $\bar{x}_{s_n} + D$.
This way, when creating the relative ordering, the neighbour is considered to be $D$ better than it actually is, which is the same as what the Simulated Annealing considers.

\subsection{$t$-Test}\label{sec:oo_ttest}
Both the OCBA and IZ procedures are designed for comparing $k$ solutions.
However, in Simulated Annealing, we only compare two solutions ($k = 2$).
Thus, we developed a third method for determining how many simulations are needed.
For this, we first look at determining a relative ordering of the solutions.
The idea here is to use a paired samples $t$-test to determine whether the two solutions can be assumed to have the same value or if they have a significant difference.
Note that we can use the paired $t$-test because of the use of CRN.
This $t$-test is repeated until a significant difference is found or $N_\text{max}$ simulations have been performed.
If a significant difference is found, we stop simulating and use the sampled averages in the Simulated Annealing to determine if the new solution is accepted.
Otherwise, the two solutions are considered the same so we perform additional simulations and repeat the test.

For this procedure, we select four parameters: the initial number of simulations $n_0$, the number of additional simulations $\Delta$, the maximum number of simulations $N_\text{max}$, and a confidence parameter $\alpha_\text{conf}$.
Note that, similar to the OCBA and IZ methods, $N_\text{max}$ denotes the maximum for the total number of simulations.
In contrast to the OCBA method, this procedure always allocates an equal number of simulations to each solution and hence a solution is simulated at most $N_\text{max} / 2$ times.
The procedure is as follows:
\begin{enumerate}
    \item Let $n$ be the number of simulations we run for a single solution and simulate each solution $n_0$ times, i.e. $n \gets n_0$.
    \item\label{item:ttest_calculate_t}
    Let $\bar{x}$ and $s^2$ be the sample mean and sample variance of the differences.
    This difference is defined as $x_{s,i} - x_{s_n,i}$ ($i = 1, ..., n$), where $x_{s,i}$ is the cost of the $i$th simulation of solution $s$.
    As defined in \Cref{sec:simulated_annealing}, solution $s$ is the current solution and $s_n$ the neighbouring solution.
    Then calculate
    \begin{equation}
        t = \frac{\bar{x}}{\sqrt{\frac{s^2}{n}}}.
    \end{equation}

    \item\label{item:ttest_decide_extra_simulations} Use the calculated value of $t$ in a two-sided $t$-test with $n - 1$ degrees of freedom to determine the $p$-value.
    Thus, \begin{equation}
        p = 2F_{n - 1}(-\abs{t}),
    \end{equation}
    where $F_{n-1}$ is the cumulative distribution of the Student's $t$-distribution with $n - 1$ degrees of freedom.
    Then, stop if $p < \alpha_\text{conf}$ or $2n \geq N_\text{max}$.

    \item Otherwise, simulate each solution an additional $\Delta$ number of times (thus $n \gets n + \Delta$) and return to Step \ref{item:ttest_calculate_t}.
\end{enumerate}

In this procedure, we use a $t$-test to determine whether the average difference is equal to zero or not.
However, as noted before, in the Simulated Annealing procedure, the new solution $s_n$ is accepted if $\hat{c}_{s_n} + D \leq \hat{c}_s$.
When $D$ is roughly the same as the cost difference between the solutions, this may result in inaccurate decisions by the Simulated Annealing, because in such cases one may not be able to accurately determine whether the difference between both solutions is greater or smaller than $D$.
Therefore, we create two different modifications to this procedure to account for the allowed difference $D$.
The first modification tests for a difference of $D$, which we later denote as $t$-test ($D$).
Our second modification adds an additional $t$-test to test whether the difference is greater than $D$ and hence we call this the double $t$-test method.

In the first modification, we test whether the difference is equal to $D$ instead of $0$.
If the test accepts this hypothesis, the Simulated Annealing might not be able to make a correct decision and hence we need more simulations.
But if the hypothesis is rejected, we know that $\bar{x}$ is either greater than $D$ or smaller than $D$.
Thus, we stop simulating as the Simulated Annealing can now make the correct decision.
To test for equality to $D$, we calculate
\begin{equation}
    t' = \frac{\bar{x} - D}{\sqrt{\frac{s^2}{n}}}
\end{equation}
and replace it with $t$ in the above procedure.

The double $t$-test method focuses on reducing the number of simulations required.
This is done by adding an additional $t$-test to the original procedure.
When two solutions are roughly equal, and hence the original procedure requires us to do additional simulations, it might be reasonable to accept $s_n$, since we accept solutions that are at most $D$ worse.
Therefore, we check whether the difference is significantly greater than $D$.
In that case, the neighbouring solution will be accepted by the Simulated Annealing.
Thus, we use the allowed difference $D$ to decide whether extra simulations are necessary.

The exact double $t$-test procedure is as follows:
\begin{enumerate}
    \item Let $n$ be the number of simulations we ran for a single solution and simulate each solution $n_0$ times, i.e. $n \gets n_0$.
    \item\label{item:ttest2_calculate_t} Let $\bar{x}$ and $s^2$ be the sample mean and sample variance of the difference, with the difference as defined before.
    Then calculate
    \begin{equation}
        t_1 = \frac{\bar{x}}{\sqrt{\frac{s^2}{n}}}
    .\end{equation}

    \item\label{item:ttest2_first_test} Use the calculated value of $t_1$ in a two-sided $t$-test with $n - 1$ degrees of freedom to determine the $p$-value. Thus, \begin{equation}
        p_1 = 2F_{n - 1}(-\abs{t_1}),
    \end{equation}
    where $F_{n-1}$ is the cumulative distribution of the Student's $t$-distribution with $n - 1$ degrees of freedom.
    Then, stop if $p_1 < \alpha_\text{conf}$ or $2n \geq N_\text{max}$.

    \item\label{item:ttest2_second_test} Perform a one-sided $t$-test to see whether the difference between the two solutions ($\bar{x}$) is significantly larger than $D$.
    Note that $D \leq 0$ and $\bar{x}$ is determined to be roughly $0$.
    Thus, we only have to test for $\bar{x} > D$, hence the one-sided $t$-test.
    For this, calculate
    \begin{equation}
        t_2 = \frac{\bar{x} - D}{\sqrt{\frac{s^2}{n}}}.
    \end{equation}
    Then,
    \begin{equation}
        p_2 = F_{n - 1}(-t_2).
    \end{equation}
    Stop if $p_2 < \alpha_\text{conf}$.
    Otherwise, continue to the next step.

    \item\label{item:ttest2_continue} Simulate each solution an additional $\Delta$ number of times (thus $n \gets n + \Delta$) and return to Step \ref{item:ttest2_calculate_t}.
\end{enumerate}

Note that if the procedure stops in Step \ref{item:ttest2_first_test}, the Simulated Annealing still needs to decide based on the sampled averages as before.
If the procedure is stopped in Step \ref{item:ttest2_second_test}, the Simulated Annealing can directly accept the new solution.

\section{Case Study: Parallel Machine Scheduling}\label{sec:cs_spmsp}
We start our computational study with a case study on the \textsc{Stochastic Parallel Machine Scheduling Problem} (SPMSP).
It is a variant of parallel machine scheduling, where a set of jobs needs to be distributed over a number of identical machines.
Then, all jobs must be completed before a given deadline.
Each job has a release date, and there are a number of precedence constraints which enforce that a job may only start once all its predecessors are finished.
The processing times of the jobs are modelled as stochastic variables.
For each processing time variable, we assume to know the type, mean and standard deviation of its probability distribution.
We construct a baseline schedule specifying the start time and machine of each job.
When the actual values of the processing times are revealed, the schedule is run according to a certain execution policy.
We consider an execution policy where jobs start as soon as possible, but never before their assigned start time.

Our objective is to find a robust schedule.
To create a robust schedule, we first need to know when a schedule is considered to be robust.
The robustness of a schedule can be characterized by two different aspects: quality robustness and solution robustness \citep{Vonder2005}.
The former indicates the stability of the objective function.
Solution robustness considers the stability of the solution itself, referring to how well the original plan can be followed.
In this problem, this would translate to the stability of the start time of each job.
Note, that a job cannot start before its scheduled start time, thus to increase solution stability one can insert buffers to prevent delay propagation.
However, this leads to a decreased probability meeting the deadline.
In the objective, we use both a quality robustness and solution robustness measure.
The quality robustness is measured by the probability of meeting the given deadline, and the solution robustness is measured by the expected fraction of jobs that start on time.

A solution to the SPMSP needs to balance these two conflicting ways to score the robustness.
For this, we developed a Simulated Annealing algorithm.
This algorithm achieves this balance by using an objective function that assigns an adaptive weight to the two robustness objectives, where the weight is based on the ratio between these two objectives in the current solution.
The neighbourhoods for generating new solutions are standard: moving a job, swapping two jobs, changing the buffer time after a job, and moving part of one job's buffer to a predecessor or successor. 
Further details are omitted for reasons of brevity.

In each iteration of the Simulated Annealing, the objective scores of the neighbour and the current solution are compared using one of the methods to compute the number of simulations (as described in \Cref{sec:number_simulations}).
Note that these methods were described with a minimization problem in mind, hence small changes are required for them to apply to a maximization problem.
For the IZ method, we need to change the $b$ and $\delta_i$ calculations to be $b = \argmax_i \bar{x}_i$ and $\delta_i = \max\{ \delta^*, \bar{x}_b - \bar{x}_i \}$.
Next to that the difference in the $t$-test methods is now defined as $x_{s_n,i} - x_{s,i}$.
These changes ensure that the solution with a bigger score is seen as the better one.

The Simulated Annealing is stopped when a certain temperature is reached.
Then, the final solution returned is evaluated with \num{10000} simulations, recording the fraction of samples meeting the deadline (quality robustness) and the fraction of jobs that can start at their planned time (solution robustness).
The average of both metrics yields the final robustness score.

The different methods for calculating the number of simulations required are compared on eight generated instances.
These instances are generated using the same method as described by \citet{Loman2025}.
Each instance is characterized by a number of jobs $j$, a number of precedence relations $r$, and the number of machines available $m$.
The processing time of each job is normally distributed around a given mean and a standard deviation which is 40\% of that mean.

The parameters used for each of the methods in \Cref{sec:number_simulations} are shown in \Cref{tab:spmsp_method_parameters}.
These parameters were found by doing parameter tuning for each of the methods.
Here, the ``Const'' methods denote methods following the simple idea of performing a constant number of simulations each iteration.
In this case, each solution is simulated $N_\text{max} / 2$ times in every iteration.
Note that the IZ and $t$-test methods also simulate each solution at most $N_\text{max} / 2$ times in every iteration.
In order to understand the effect of CRN, we included a ``Const'' method with CRN disabled.
In the results, the IZ and $t$-test methods contain either the ($0$) or ($D$) suffix to denote if they consider a relative ordering or a relative ordering with respect to the allowed difference $D$.

\begin{table}[ht!]
    \centering
    \caption{Parameters used for the different methods for calculating the required number of simulations for the SPMSP.}
    \label{tab:spmsp_method_parameters}
    \begin{tabular}{lrrrr}
        \toprule
        Method          & $n_0$ & $\Delta$ & $N_\text{max}$ & $\alpha_\text{conf}$ \\
        \midrule
        Const           &       &          &            100 &          \\
        Const           &       &          &            200 &          \\
        Const           &       &          &            400 &          \\
        Const (no CRN)  &       &          &            400 &          \\
        OCBA            &    80 &       10 &            400 &          \\
        IZ              &    80 &       10 &            400 &    $0.2$ \\
        TTest           &    80 &       20 &            400 &    $0.2$ \\
        Double-TTest    &    80 &       20 &            400 &    $0.2$ \\
        \bottomrule
    \end{tabular}
\end{table}

\begin{figure}[ht!]
    \centering
    \includegraphics[width=0.95\textwidth]{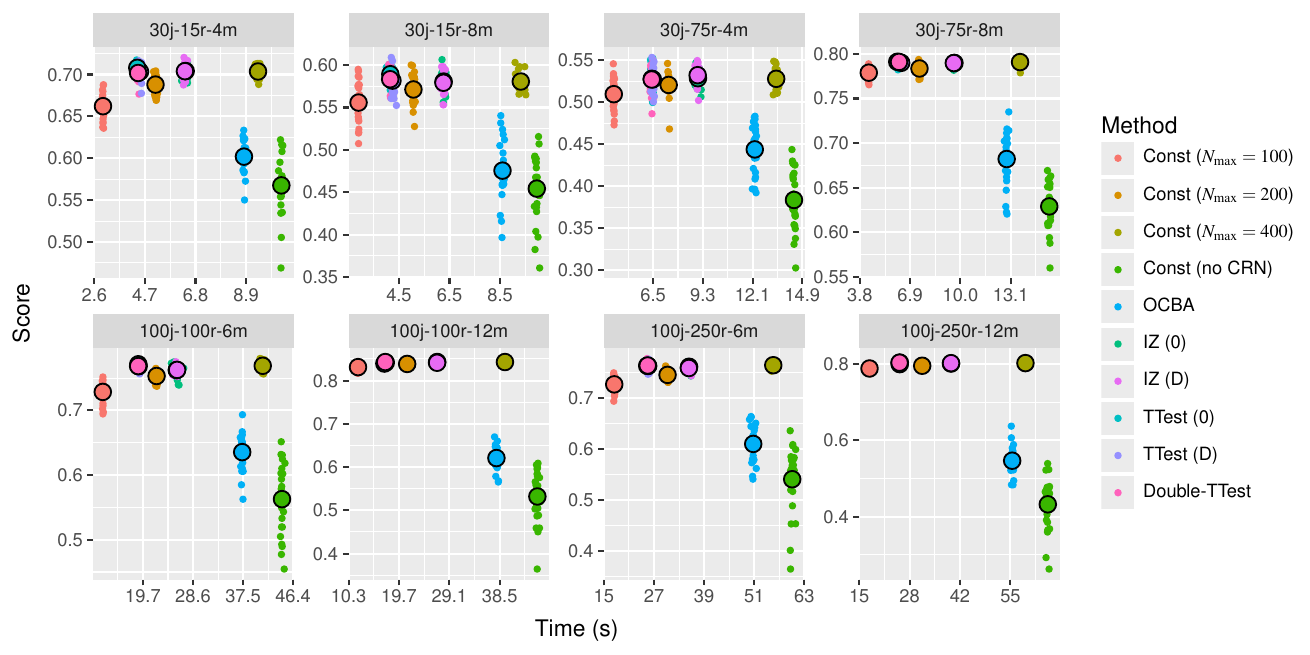}
    \caption{Comparison of the robustness score and runtime of various Simulated Annealing runs for different methods for calculating the required number of simulations.}
    \label{fig:spmsp_score_time_no_crn}
\end{figure}

First, we look at the performance of each method.
For each method, 25 independent runs of Simulated Annealing are performed.
The results of these 25 runs are plotted in \Cref{fig:spmsp_score_time_no_crn}, showing the runtime and achieved robustness score of each of these runs.
Note that this robustness score is the average between the probability of meeting the given deadline, and the expected fraction of jobs that start on time.
The bigger dots in this image show the average runtime and average score over the 25 runs of each method.
In here, we notice that the methods that do not employ CRN (OCBA and Const without CRN) score significantly worse than the other methods.
Not using CRN results in less fair comparisons and the results show that this leads to incorrect decisions.
These incorrect decisions lead to both a worse objective value and a higher variance of this objective.
This indicates the need to employ CRN when comparing two solutions.

\begin{figure}[ht!]
    \centering
    \includegraphics[width=0.95\textwidth]{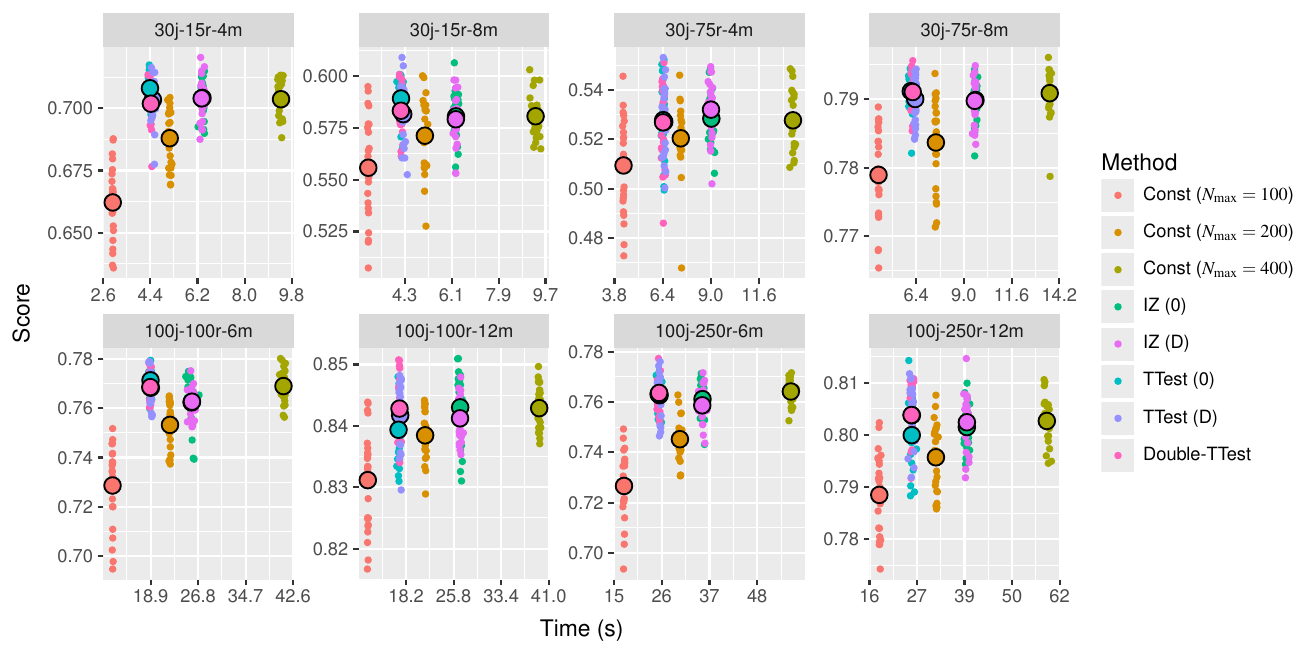}
    \caption{Comparison of the robustness score and runtime of various Simulated Annealing runs for different methods for calculating the required number of simulations. Here, all methods that do not employ CRN are excluded.}
    \label{fig:spmsp_score_time_crn}
\end{figure}

\begin{figure}[ht!]
    \centering
    \includegraphics[width=0.95\textwidth]{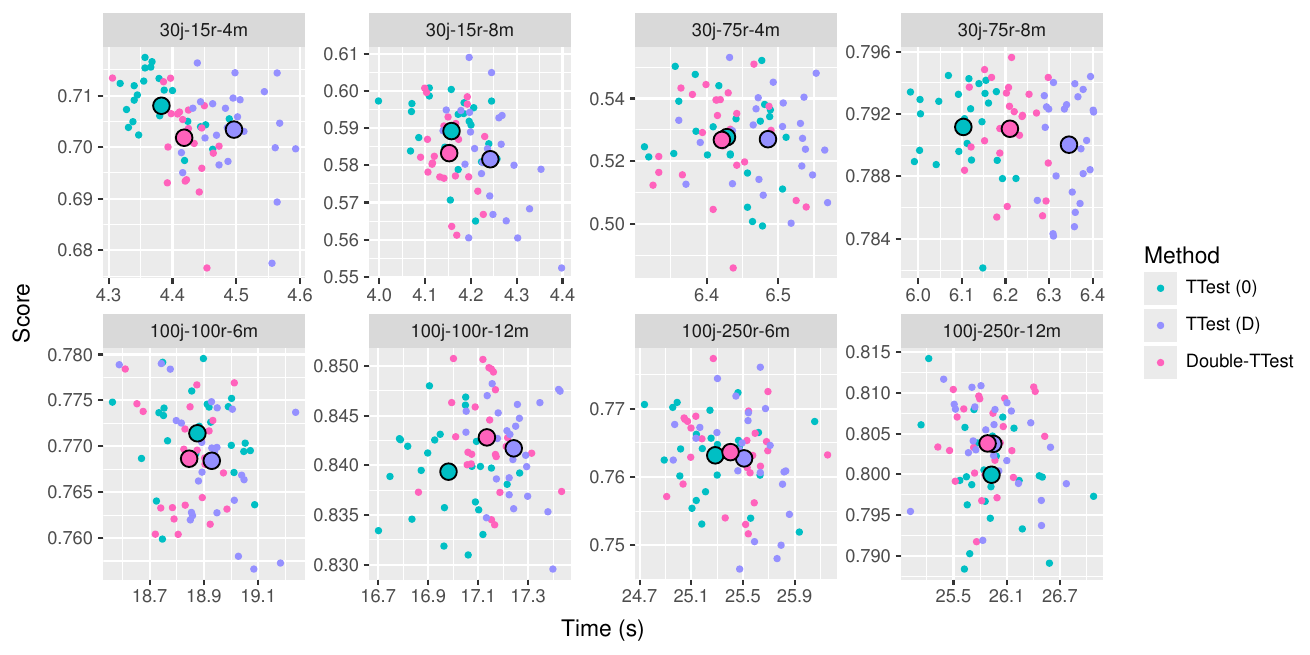}
    \caption{Comparison of the robustness score and runtime of various Simulated Annealing runs for the different $t$-test methods.}
    \label{fig:spmsp_score_time_ttest}
\end{figure}

To get a better picture of the methods that do employ CRN, we compare them in \Cref{fig:spmsp_score_time_crn}.
In terms of solution quality, the $t$-test methods and the IZ method outperform the ``Const'' methods.
For the IZ method, there seem to be minimal differences between using the allowed difference information or not.
In some instances ``IZ ($D$)'' is a bit better, while in others it is not.
The runtime of the IZ methods show an improvement compared to using an equal distribution with $N_\text{max} = 400$, while getting similar scores.
The $t$-test methods themselves are further highlighted in \Cref{fig:spmsp_score_time_ttest}.
These methods show very similar performances, where it is not clear if there are significant differences between them.
We do note that ``TTest ($D$)'' seems to have the slowest runtime on average, but there does not seem to be a method that is clearly dominating the others.

\begin{figure}[ht!]
    \centering
    \includegraphics[width=0.80\textwidth]{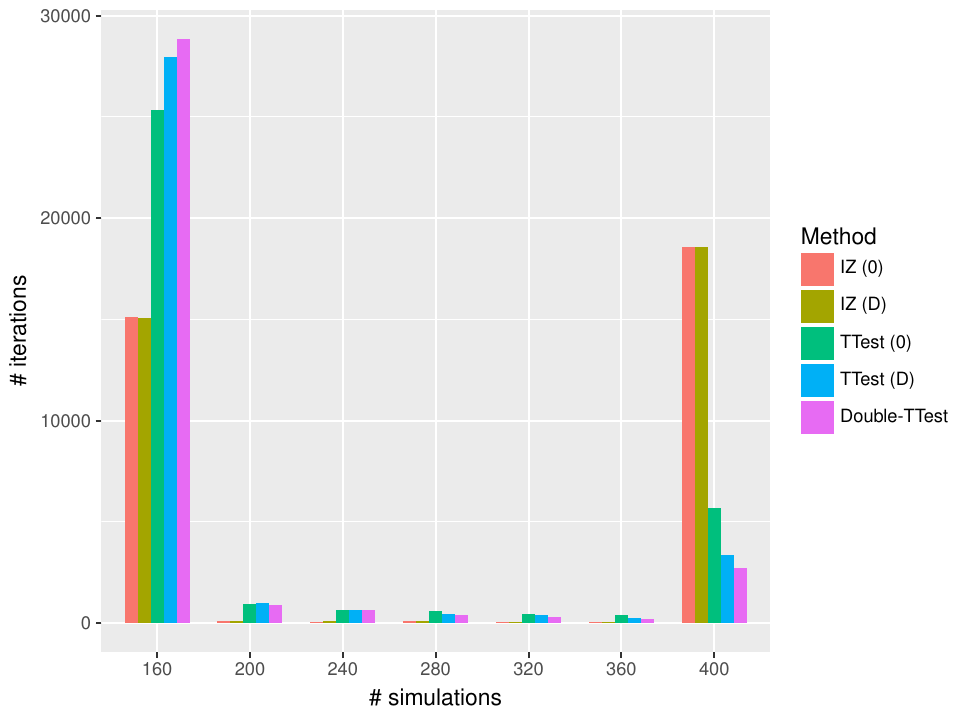}
    \caption{Histogram of the number of simulations done in a single Simulated Annealing iteration.}
    \label{fig:spmsp_sim_dist}
\end{figure}

To get a better understanding of the runtime performance of the IZ and $t$-test methods, we investigate the number of simulations done in a single iteration of Simulated Annealing.
A histogram of this is shown in \Cref{fig:spmsp_sim_dist}.
Here we again see no clear difference between ``IZ ($0$)'' and ``IZ ($D$)'', with both using a similar number of simulations each iteration, often using the maximum number allowed.
This explains the difference in runtime between IZ and the $t$-test methods.
Between the $t$-test methods, we notice that using the allowed difference information has a big impact on the number of simulations performed.
The ``TTest ($D$)'' method uses fewer simulations than ``TTest ($0$)'' and the ``Double-TTest'' uses the fewest overall.
However, as can be seen in \Cref{fig:spmsp_score_time_ttest}, this does not translate to a clear overall speed advantage.

Finally, we look at the convergence of each method for calculating the number of required simulations.
For this, we track the score of the best found solution in each iteration of the Simulated Annealing.
Each time $s_\text{best}$ is updated, we simulate the new $s_\text{best}$ \num{1000} times and log the average score.

\begin{figure}[ht!]
    \centering
    \includegraphics[width=0.80\textwidth]{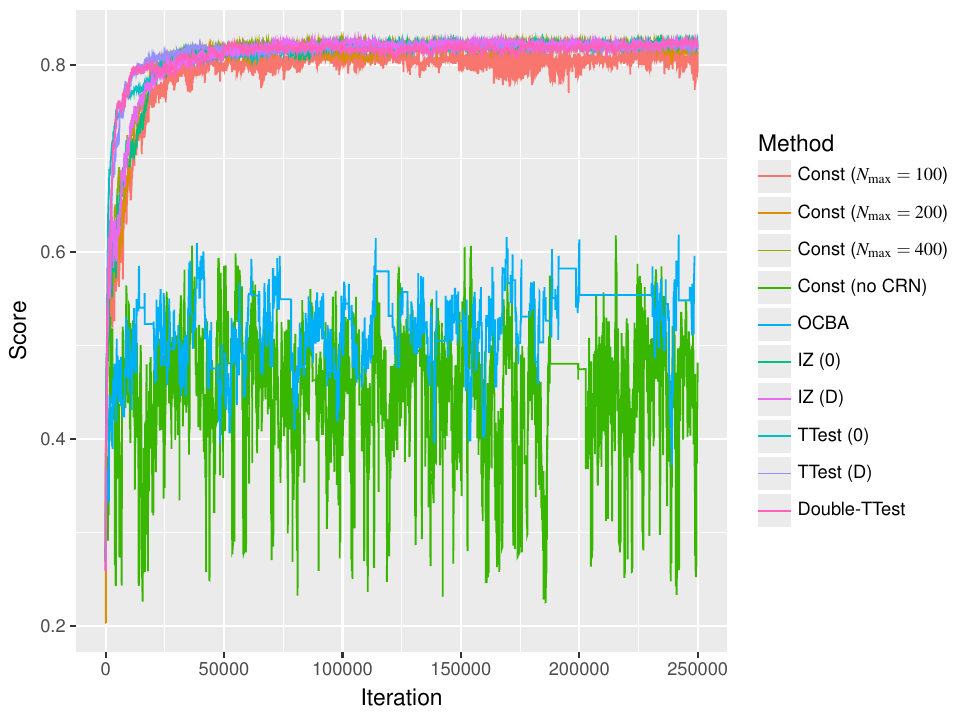}
    \caption{Best solution over time for each comparison method.}
    \label{fig:spmsp_convergence}
\end{figure}

The first \num{250000} iterations of the Simulated Annealing are shown in \Cref{fig:spmsp_convergence}.
Here we show the results of the ``100j-250r-12m'' instance.
In this figure, we see that the methods employing CRN show convergence, while the methods not employing CRN suffer from poor convergence.
Moreover, this figure reveals that the best solution is sometimes worse than the previous one.
This indicates that the method used for determining the new best solution does not always make a correct decision.
Note that when the algorithm converges, the best solution will only improve in very small steps and hence it is more difficult to distinguish the previous and next best solutions.
Still, the differences in the graph are small, especially if $N_\text{max}$ is larger.

\section{Case Study: Stochastic E-VSP}\label{sec:cs_evsp}
In the second case study, we study the \textsc{Vehicle Scheduling Problem} for public bus transport.
We specifically look at the \textsc{Electric Vehicle Scheduling Problem} (E-VSP).
Here, we are given a set of trips from the timetable, with their departure and arrival locations, start times, and driving times.
The goal is then to schedule electric buses, i.e. assign a sequence of trips to each bus, such that every trip is driven and the cost is minimized.
Since electric buses cannot drive the entire day without charging, the route of a single vehicle has to contain charging.
An overview of this problem is provided in a recent survey by \citet{Perumal2021}.

In this case study, we follow the work of \citet{Bosch2021} and \citet{Bruin2023}.
\Citet{Bosch2021} propose a Simulated Annealing algorithm to solve the problem where they consider costs for battery degradation as well as non-linear charging times for buses.
\Citet{Bruin2023} considered the stochastic variant of this problem, where they considered both stochastic driving times and stochastic energy consumption.
Distributions for these stochastic variables were derived from historical driving times and industry insights.
They add a penalty cost for trips that are started late.
Thus, the objective not only minimizes total cost, but also expected lateness.
Then, they use Discrete-Event Simulation to estimate the objective function.
These estimations were then used in the Simulated Annealing approach of \citet{Bosch2021}.

Here, we compare the different methods for computing the required number of simulations on the same algorithm and instances as used by \citet{Bruin2023}.
These instances are from 4 areas as served by Qbuzz, a major bus operator in the Netherlands.
The sizes and properties of these instances are shown in \Cref{tab:evsp_instance_overview}.

\begin{table}[ht!]
    \centering
    \caption{Overview of the used datasets and their parameters.}
    \label{tab:evsp_instance_overview}
    \begin{tabular}{lrrrr}
        \toprule
        Dataset          & \#Trips & \#Lines & Battery Capacity (kWh) \\
        \midrule
        \texttt{dmg}     &     631 &       8 &                    232 \\
        \texttt{gn345}   &     463 &       3 &                    184 \\
        \texttt{qlink}   &     590 &       3 &                    160 \\
        \texttt{zst}     &     317 &       2 &                    232 \\
        \bottomrule
    \end{tabular}
\end{table}

The parameters used for each of the methods in \Cref{sec:number_simulations} are shown in \Cref{tab:evsp_method_parameters}.
As in the previous case study, note that for the ``Const'' methods, each instance gets simulated $N_\text{max} / 2$ times in every iteration of the Simulated Annealing.
Furthermore, in order to understand the effect of CRN, we included a ``Const'' method with CRN disabled.
In the results, the IZ and $t$-test methods contain either the ($0$) or ($D$) suffix to denote if they consider a relative ordering or a relative ordering with respect to the allowed difference $D$.

\begin{table}[ht!]
    \centering
    \caption{Parameters used for the different methods for calculating the required number of simulations.}
    \label{tab:evsp_method_parameters}
    \begin{tabular}{lrrrr}
        \toprule
        Method          & $n_0$ & $\Delta$ & $N_\text{max}$ & $\alpha$ \\
        \midrule
        Const           &       &          &             20 &          \\
        Const           &       &          &             50 &          \\
        Const           &       &          &            200 &          \\
        Const (no CRN)  &       &          &            200 &          \\
        OCBA            &    20 &        5 &            200 &          \\
        IZ              &    10 &       10 &            200 &    $0.1$ \\
        TTest           &    10 &       10 &            200 &    $0.2$ \\
        Double-TTest    &    10 &       10 &            200 &    $0.2$ \\
        \bottomrule
    \end{tabular}
\end{table}

\begin{figure}
    \centering
    \includegraphics[width=0.95\textwidth]{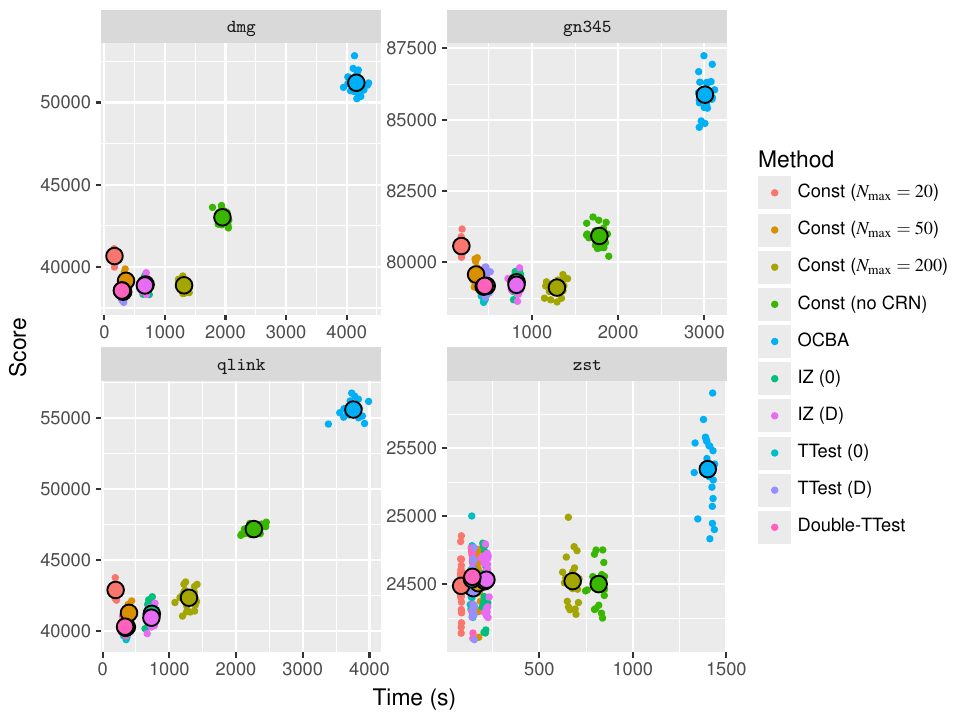}
    \caption{Comparison of the score and runtime of various Simulated Annealing runs for different methods of calculating the required number of simulations.}
    \label{fig:evsp_score_time_no_crn}
\end{figure}

As in \Cref{sec:cs_spmsp}, we first compare the runtime performance and the solution scores.
For this, we perform 25 runs of Simulated Annealing and compare these runs in \Cref{fig:evsp_score_time_no_crn}.
Note that the E-VSP problem is a minimization problem, thus lower scores are better.
Furthermore, the bigger dots, again, show the average runtime and average score over the 25 runs of each method.
Similar to \Cref{fig:spmsp_score_time_no_crn}, we see a clear performance difference between the methods with CRN and the methods without.
This further shows the need to employ CRN.

\begin{figure}
    \centering
    \includegraphics[width=0.95\textwidth]{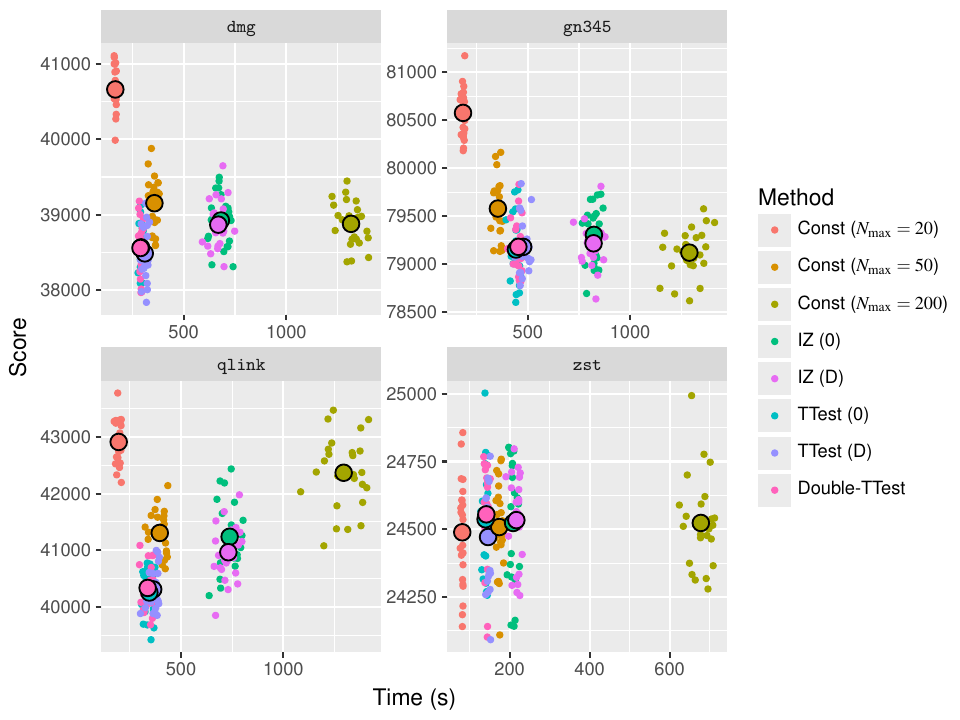}
    \caption{Comparison of the score and runtime of various Simulated Annealing runs for different methods of calculating the required number of simulations. Here, all methods that do not employ CRN are excluded.}
    \label{fig:evsp_score_time_crn}
\end{figure}

\begin{figure}[ht!]
    \centering
    \includegraphics[width=0.95\textwidth]{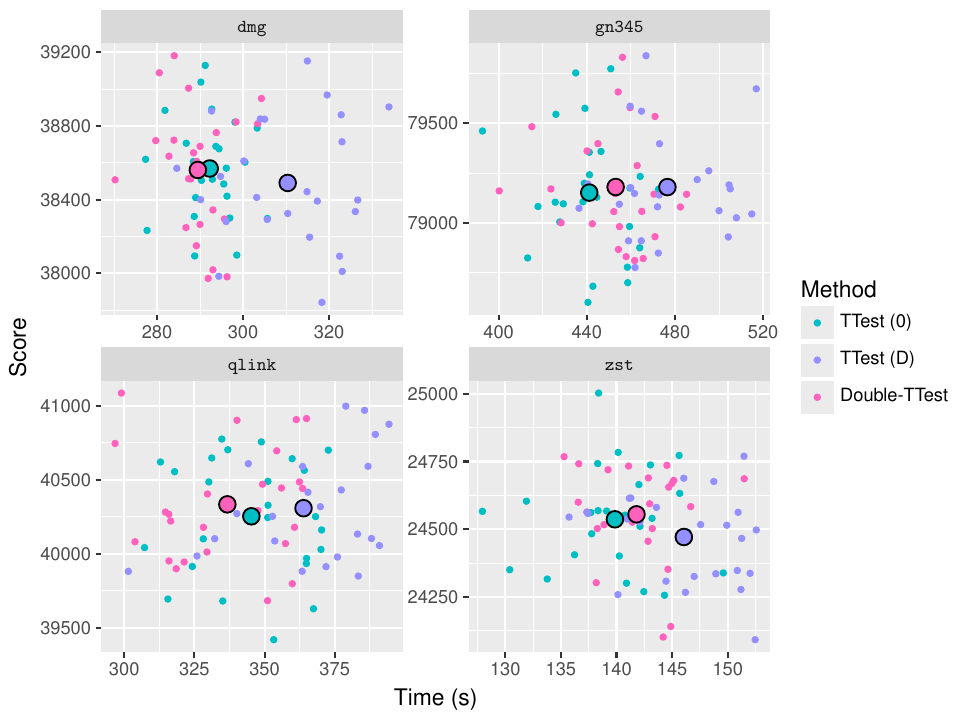}
    \caption{Comparison of the score and runtime of various Simulated Annealing runs for the different $t$-test methods.}
    \label{fig:evsp_score_time_ttest}
\end{figure}

In \Cref{fig:evsp_score_time_crn}, we compare only the methods that do employ CRN.
The \texttt{zst} instance, which is the smallest, is the easiest to solve and all methods return similar solution outcomes, albeit with different runtimes.
For the other instances a similar picture is visible as in \Cref{fig:spmsp_score_time_crn}.
However, in this case there is a clearer distinction between ``IZ ($0$)'' and ``IZ ($D$)'', with  the latter being faster but on average a bit worse.
The $t$-test methods show the best performance when it comes to score and runtime.
A comparison between the $t$-test methods is made in \Cref{fig:evsp_score_time_ttest}.
Similar to \Cref{fig:spmsp_score_time_ttest}, we note that ``TTest ($D$)'' seems to be the slowest of the three, while all three methods produce similar solutions cost-wise.

\begin{figure}
    \centering
    \includegraphics[width=0.80\textwidth]{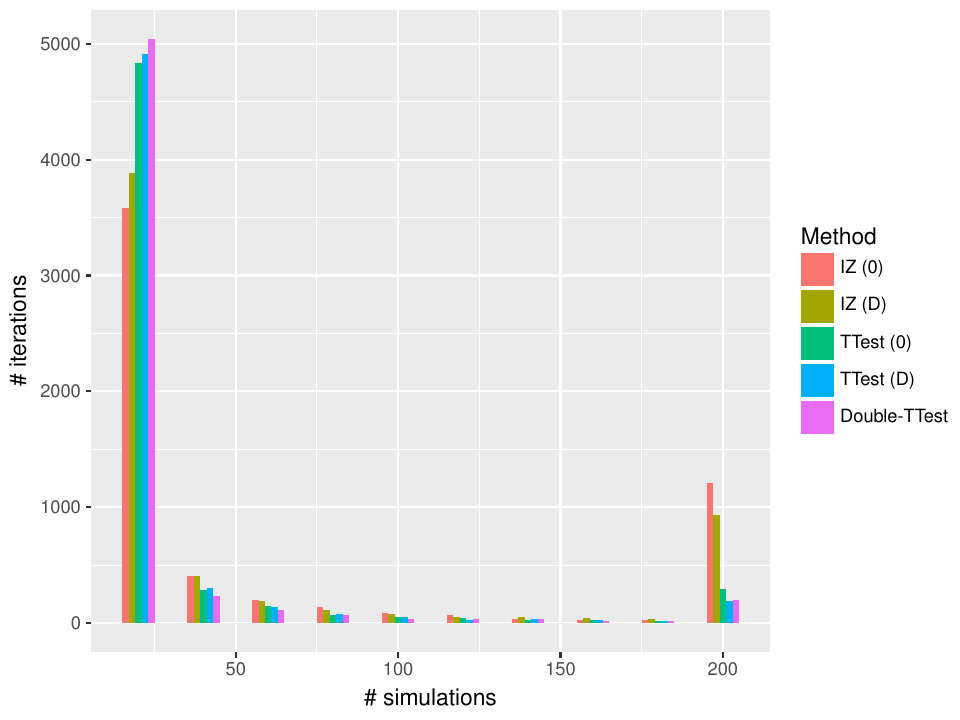}
    \caption{Histogram of the number of simulations done in an iteration of Simulated Annealing.}
    \label{fig:evsp_sim_dist}
\end{figure}

Next, we show the number of simulations performed in an iteration of local search in \Cref{fig:evsp_sim_dist}.
Contrary to \Cref{fig:spmsp_sim_dist}, there is a bigger difference between the ``IZ ($0$)'' and ``IZ ($D$)'' methods, with the latter using fewer simulations.
The $t$-test methods, however, perform similar to the SPMSP case, with the Double-TTest using the fewest simulations overall.

\begin{figure}
    \centering
    \includegraphics[width=0.80\textwidth]{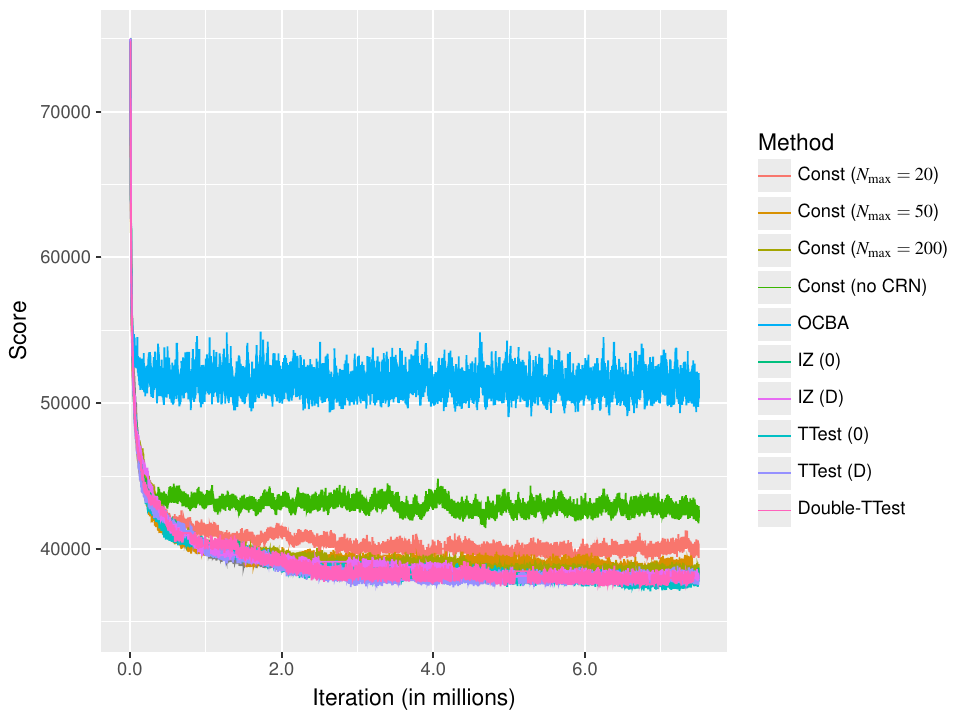}
    \caption{Best solution over time for each comparison method.}
    \label{fig:evsp_convergence}
\end{figure}

Lastly, we look at the convergence of each method for calculating the number of required simulations.
We do this in a similar manner as in \Cref{sec:cs_spmsp}, where we keep track of the score of the current best solution.
Each time the best solution is updated, we simulate the new best solution \num{100} times and log the score.
The first $7.5$ million iterations of the Simulated Annealing are shown in \Cref{fig:evsp_convergence}.
These are the results of the \texttt{dmg} instance.
The results seem less extreme than in the SPMSP case study (\Cref{fig:spmsp_convergence}), although the non-CRN methods still perform quite bad as they show a large variance in their best solution scores.
The results of the other methods are similar to \Cref{fig:spmsp_convergence}, where the errors of wrongly accepting a new best solution appear to become smaller for methods with a higher $N_\text{max}$.

\section{Conclusion}\label{sec:conclusion}
In this paper, we investigate the use of simulations inside a local search algorithm, with a specific focus on Simulated Annealing.
An important aspect here is the number of simulations performed in a single iteration of the local search.
We propose a new method based on $t$-tests to reduce the number of simulations and compare it to existing methods.
These methods are compared on the Stochastic Parallel Machine Scheduling Problem and the stochastic Electric Vehicle Scheduling Problem.

In both our case studies, we see a similar pattern of results regarding runtime and solution quality.
Our experiments show the benefit of employing Common Random Numbers (CRN), as not using CRN results in longer runtimes, worse solutions, and bad convergence for both methods we tested that do not employ CRN.
This behaviour seems intuitive, since employing CRN will create a fairer comparison, especially when dealing with a relatively low number of simulations.
Our newly proposed $t$-test methods outperform the use of Indifference Zones, showing better runtimes and creating solutions with equal or better quality.
This highlights the benefit of using one of the $t$-test methods to determine the number of simulations required.

The $t$-test methods themselves perform very similarly to each other.
In both case studies, they show a similar runtime and solution outcome.
Here, we notice that the use of the ``Double-TTest'' method reduces the number of simulations performed, however the runtime cost of doing extra $t$-tests seems similar to the runtime of these simulations.
In other words, the time saved by not performing these simulations is compensated by performing extra $t$-tests.
This results in the ``Double-TTest'' method having a similar runtime as the other two $t$-test methods.
However, in other problems this behaviour could be different as it depends on how efficient the simulation can be implemented.

Next to that, we adapted the IZ and $t$-test methods to take into account that Simulated Annealing can also accept worse solutions.
The results for these adapted methods were not always intuitive.
For the IZ method, it did not seem to make much of a difference, while for the $t$-test method, it was often the slowest, although not significantly slower than the other $t$-test methods.
It remains hard to explain why exactly this is the case.

In our setup, we require the use of simulations every iteration of the local search.
However, there might be ways to integrate with ad hoc approximations such that simulations may not be required every iteration.
This raises the question of how and if, for example, robustness measures can be included to reduce total simulations while keeping solution quality the same.
One idea for this would be to replace the step of performing $n_0$ simulations with the evaluation of a robustness measure.
The question then remains how the information from this measure could be used to determine whether extra simulations are required.
Answering these questions is an interesting topic for further research.

\backmatter

\bibliography{papers}

\end{document}